\documentclass[reqno,12pt]{amsart}

\usepackage{amssymb}
\usepackage{amsthm}

\numberwithin{equation}{section}

\newtheorem{theorem}{Theorem}[section]

\newtheorem{lemma}[theorem]{Lemma}
\newtheorem{corollary}[theorem]{Corollary}

\theoremstyle{definition}

\newtheorem{remark}[theorem]{Remark}
\newtheorem{remarks}[theorem]{Remarks}
\newtheorem{example}[theorem]{Example}

\newcommand{\R}{\mathbb{R}}

\newcommand{\N}{\mathbb{N}}

\newcommand{\M}{\mathbb{M}}
\newcommand{\Sym}{\mathbb{S}}
\newcommand{\Ort}{\mathbb{O}}

\newcommand{\SO}{\mathbb{SO}}
\newcommand{\GL}{\mathbb{GL}}

\newcommand{\cT}{\mathcal{T}}
\newcommand{\cNT}{\mathcal{NT}}
\newcommand{\cNS}{\mathcal{NS}}
\newcommand{\cS}{\mathcal{S}}
\newcommand{\cE}{\mathcal{E}}
\newcommand{\cH}{\mathcal{H}}

\newcommand{\cMh}{\mathcal{M_H}}

\newcommand{\diag}[1]{\mathrm{diag}\left(#1\right)}
\newcommand{\sgn}[1]{\mathrm{sgn}#1}
\newcommand{\tr}[1]{\mathrm{tr}(#1)}
\newcommand{\ind}[2]{\mathrm{i}\left(#1,#2\right)}
\newcommand{\dist}[2]{\mathrm{dist}\left(#1,#2\right)}
\newcommand{\spect}[1]{\sigma\left(#1\right)}
\newcommand{\pspect}[1]{\sigma_{+}\left(#1\right)}
\newcommand{\nspect}[1]{\sigma_{-}\left(#1\right)}
\newcommand{\eigsp}[2]{V_{#1}\left(#2\right)}
\newcommand{\mult}[1]{\mu(#1)}
\newcommand{\morse}[1]{\mathrm{m}^{-}\left(#1\right)}
\newcommand{\halfline}{(0,+\infty)}
\newcommand{\abs}[1]{\left\vert#1\right\vert}
\newcommand{\norm}[1]{\left\Vert#1\right\Vert}
\newcommand{\set}[1]{\left\{#1\right\}}
\newcommand{\cl}[1]{\mathrm{cl}\left(#1\right)}
\newcommand{\cond}{\;\vert\;\;\;}
\newcommand{\hilb}{H^{1}_{2\pi}}

\begin{document}


\title[Degenerate Bifurcation Points]{Degenerate Bifurcation Points
of Periodic Solutions of Autonomous Hamiltonian Systems}


\author[W. Radzki]{Wiktor Radzki$^{\dag}$}


\address{$^{\dag}$Faculty of Mathematics and Computer Science,
  Nicolaus Copernicus University,
  ul. Chopina $12 \slash 18$,
  87-100 Toru\'{n},
  Poland}
\email{radzki@mat.uni.torun.pl}


\author[S. Rybicki]{S{\l}awomir Rybicki$^{\ddag}$}


\address{$^{\ddag}$Faculty of Mathematics and Computer Science,
  Nicolaus Copernicus University,
  ul. Chopina $12 \slash 18$,
  87-100 Toru\'{n},
  Poland}
\email{rybicki@mat.uni.torun.pl}
\thanks{$^{\ddag}$Research supported by   the
SCSR under grant number  5 PO3A 026 20.}


\keywords{Hamiltonian system, Periodic solution, Bifurcation,
Emanation, Branching point, Bifurcation index, Topological degree
for \mbox{$\mathbb{SO}(2)$-equi}variant gradient maps}


\subjclass{Primary: 34C23, 34C25; Secondary: 58F14, 70H05.}


\date{March 28, 2004}


\begin{abstract}
We study connected branches of non-constant
\mbox{$2\pi$-pe}riodic solutions of the Hamilton equation
\begin{displaymath}
\dot{x}(t)=\lambda J\nabla H(x(t)),
\end{displaymath}
where $\lambda\in\halfline,$ $H\in C^2(\R^n\times\R^n,\R)$ and
$
\displaystyle \nabla^2H(x_0)=
\left[
\begin{array}{cc}
A&0\\
0&B
\end{array}
\right]
$
for $x_0\in\nabla H^{-1}(0).$ The Hessian $\nabla^2H(x_0)$ can be
singular. We formulate sufficient conditions for the existence
of such branches bifurcating from given $(x_0,\lambda_0).$
As a consequence we prove theorems concerning the existence
of connected branches of arbitrary periodic nonstationary
trajectories of the Hamiltonian system $\dot{x}(t)=J\nabla H(x(t))$
emanating from $x_0.$ We describe also minimal periods
of trajectories near $x_0.$
\end{abstract}


\maketitle


\section{Introduction}


Consider the autonomous Hamiltonian system
\begin{equation}{\label{ham}}
\dot{x}(t)=J\nabla H(x(t)),
\end{equation}
where $H\in C^2(\R^n\times\R^n,\R)$ and $J$ is the standard
\mbox{$2n$-di}mensional simplectic matrix.
The problem of finding periodic solutions of \eqref{ham} is
equivalent to the problem of finding solutions of the family
\begin{equation}{\label{parham}}
\left\{
\begin{array}{l}
\dot{x}(t)=\lambda J\nabla H(x(t))\\
x(0)=x(2\pi),
\end{array}
\right.
\end{equation}
with $\lambda\in\halfline.$
Let $(x_0,\lambda_0)\in \nabla H^{-1}(0)\times\halfline,$
where $\nabla H^{-1}(0)\equiv (\nabla H)^{-1}(\set{0}).$
Having a connected branch of nontrivial solutions
of \eqref{parham} bifurcating (in a suitable space) from
$(x_0,\lambda_0)$ we can find the corresponding
connected branch of nonstationary periodic trajectories
of \eqref{ham} emanating from $x_0$ with periods tending to
$2\pi\lambda_0$ at $x_0.$ Our aim is to study
such connected branches of bifurcations and emanations when
the Hessian of $H$ at $x_0$ has the block-diagonal form
\begin{equation}{\label{hessian}}
\nabla^2H(x_0)=
\left[
\begin{array}{cc}
A&0\\
0&B
\end{array}
\right],
\end{equation}
where $A$ and $B$ are real symmetric \mbox{$(n\times n)$-ma}trices.
The critical point $x_0$ can be degenerate, i.e.
$\nabla^2H(x_0)$ can be singular. However, we assume
that $x_0$ is isolated in $\nabla H^{-1}(0)$ and
the Brouwer degree of $\nabla H$ around $x_0$ is nonzero.


Condition \eqref{hessian} is satisfied, for example, if
\begin{equation}{\label{form}}
H(x)=H(y,z)=
\frac{1}{2}\langle M^{-1}y,y\rangle +V(z),
\end{equation}
where $y,z\in\R^n,$ $V\in C^2(\R^n,\R)$ and $M$ is nonsingular
real symmetric \mbox{$(n\times n)$-ma}trix.
In such a case Eq. \eqref{ham} is equivalent
to the Newton equation
\begin{equation}{\label{newton}}
M\ddot{z}(t)=-\nabla V(z(t)).
\end{equation}


Two basic results concerning bifurcations of periodic solutions
from a nondegenerate stationary point of Hamiltonian system are
due to Liapunov and Berger. If $J\nabla ^2H(x_0)$ is nonsingular
and has two purely imaginary eigenvalues $\pm i\beta$ of
multiplicity 1 then Liapunov center theorem ensures the existence
of a one-parameter family of nonstationary periodic solutions of
\eqref{ham} emanating from nondegenerate $x_0\in\nabla H^{-1}(0)$
(see \cite{L}).
Berger \cite{Br1,Br2} proved the existence of a sequence
of nonstationary periodic solutions convergent to the
nondegenerate stationary solution of \eqref{newton} for $M=I,$
without any assumptions on multiplicity of eigenvalues of the
Hessian of $V$ (see also \cite{B,T}). The above results were
generalized in \cite{DR,S,Z} to the case of Hamiltonian systems
with degenerate stationary points. Authors of \cite{S,Z} used Morse
theory and they obtained, similarly as Berger, sequences of
periodic solutions. In \cite{DR} connected branches of
nontrivial solutions of \eqref{parham} were obtained by using
the topological degree theory for \mbox{$\mathbb{SO}(2)$-equi}variant
gradient maps (see \cite{D}) and the results from \cite{RA}.
Global bifurcation theorems of this type can be also
found in \cite{GM,RY}. In the present paper we apply results from
\cite{DR} to prove the existence of connected branches of
bifurcations and emanations in degenerate case under assumptions
written in terms of the topological degree of $\nabla H$ around
$x_0$ and eigenvalues of $\nabla^2H(x_0).$ In particular, we
generalize results from \cite{Rd}, proved for $H$ satisfying
\eqref{hessian} with $A=I.$


\section{Preliminaries}


In this section we set up notation and summarize without proofs
the relevant material on bifurcation theory for Hamiltonian systems.


The number of elements of any finite set $X$ will be denoted
as $\sharp X.$
Write $\M(n,\R)$ for the set of real
\mbox{$(n\times n)$-ma}trices.
Let $I_n\in\M(n,\R)$ stand for the identity matrix. Define
$J_n\in\M(2n,\R)$ as
\begin{displaymath}
J_n=
\left[
\begin{array}{cc}
0&-I_n\\
I_n&0
\end{array}
\right].
\end{displaymath}
Usually we abbreviate $I_n$ and $J_n$ to $I$ and $J$.
Write $\Sym(n,\R),$ $\Ort(n,\R)$ and $\GL(n,\R)$ for the subsets of
$\M(n,\R)$ consisting of symmetric, orthogonal and nonsingular
matrices, respectively. If $\alpha_1,\ldots,\alpha_n\in\R$ then
$\diag{\alpha_1,\ldots,\alpha_n}$ denotes the diagonal matrix with
$\alpha_1,\ldots,\alpha_n$ on the main diagonal.


Let $\spect{A}$ be the spectrum of $A\in\M(n,\R).$
Obviously, if $A\in\Sym(n,\R)$ then $\spect{A}\subset\R.$ Symbols
$\pspect{A}$ and $\nspect{A}$ denote the sets of strictly positive
and strictly negative real eigenvalues of $A,$ respectively.
Let $\mult{\alpha}$ be the multiplicity of the eigenvalue
$\alpha\in\spect{A}$ and write $V_A(\alpha)$ for the eigenspace
corresponding to $\alpha.$
For symmetric $A$ we define the
\emph{Morse index $\morse{A}$} of $A$ by the formula
\begin{displaymath}
\morse{A}=\sum_{\alpha\in\nspect{A}}\mult{\alpha}.
\end{displaymath}


If $B\in\Sym(n,\R)$ is nonnegative definite then there exists
a nonnegative definite $C\in\Sym(n,\R)$ such that $C^2=B.$ We
denote $C$ by $\sqrt{B}.$ For any $A\in\Sym(n,\R)$ write
$\abs{A}=\sqrt{A^2}.$ If $A$ is nonnegative or nonpositive
definite then we write
\begin{displaymath}
\sgn(A)=
\left\{
\begin{array}{ll}
1;&\nspect{A}=\emptyset,\;A\neq 0,\\
0;&A=0,\\
-1;&\pspect{A}=\emptyset,\;A\neq 0.
\end{array}
\right.
\end{displaymath}
In such a case $A=\sgn(A)\abs{A}.$


For any $K,L\in\Sym(2n,\R)$ and $j\in\N$ set
\begin{displaymath}
Q_j(K)=
\left[
\begin{array}{cc}
-K&jJ^t\\
jJ&-K
\end{array}
\right],
\end{displaymath}
\begin{displaymath}
\Lambda_j(L)=
\set{\lambda\in\halfline\cond \det{Q_j(\lambda L)=0}},
\end{displaymath}
\begin{displaymath}
\Lambda(L)=\bigcup_{j\in\N}\Lambda_j(L).
\end{displaymath}
Obviously, $Q_j(K)\in\Sym(4n,\R)$ and since $\Lambda_j(L)$ is
a subset of roots of a polynomial, it is finite.
Moreover, observe that the following lemma holds true.


\begin{lemma}{\label{finite}}
For every $j\in\N$ we have
\begin{displaymath}
\Lambda_j(L)=
j\Lambda_1(L)\equiv\set{j\alpha\cond \alpha\in\Lambda_1(L)}.
\end{displaymath}
In particular,
for every $a,b\in\R,$ $a<b,$ the set $\Lambda(L)\cap[a,b]$ is finite.
\end{lemma}


In what follows we assume that $H\in C^2(\R^n\times\R^n,\R).$
We call $\cT(H)=\nabla H^{-1}(0)\times\halfline$ the set of
\emph{trivial solutions} of \eqref{parham}. The set $\cNT(H)$ of
\emph{nontrivial solutions} of \eqref{parham} consists of those
solutions $(x,\lambda)$ of \eqref{parham} that do not belong to
$\cT(H).$ We will consider $\cT(H)$ and $\cNT(H)$ as subsets of
$\hilb\times\halfline,$ see \cite{MW} for the definition of
the Hilbert space $\hilb\equiv H^1([0,2\pi],\R^{2n}).$
(Recall that $x(0)=x(2\pi)$
for every $x\in\hilb.$) The trivial solution
$(x_0,\lambda_0)\in\cT(H)$ is said to be a \emph{bifurcation point
of nontrivial solutions} of \eqref{parham} if it is a cluster
point of $\cNT(H).$ We say that a set $C\subset\cNT(H)$
\emph{bifurcates} from $(x_0,\lambda_0)$ if
$(x_0,\lambda_0)\in\cl{C}.$
Denote by $C(x_0,\lambda_0)$ the
connected component of $\cl{\cNT(H)}$ containing the bifurcation
point $(x_0,\lambda_0).$ If
$C(x_0,\lambda_0)\neq\set{(x_0,\lambda_0)}$ then $(x_0,\lambda_0)$
is called a \emph{branching point of nontrivial solutions}
of \eqref{parham}.
It was proved in \cite{DR} that any point
$(x,0),$ $x\in\hilb,$ cannot be a cluster point of $\cNT(H)$ in
$\hilb\times\R.$

The following theorem (see \cite{DR}) gives a necessary
condition for $(x_0,\lambda_0)$ to be a bifurcation point.


\begin{theorem}{\label{nec}}
Let $x_0$ be isolated in $\nabla H^{-1}(0).$
If $(x_0,\lambda_0)\in\cT(H)$ is a bifurcation
point of nontrivial solutions of \eqref{parham} then
$\lambda_0\in\Lambda(\nabla^2H(x_0)).$
\end{theorem}


Fix isolated $x_0\in\nabla H^{-1}(0),$
$\lambda_0\in\Lambda(\nabla^2H(x_0))$ and choose $\varepsilon>0$ such
that $\lambda_0-\varepsilon>0$ and
$\Lambda(\nabla^2H(x_0))
\cap [\lambda_0-\varepsilon,\lambda_0+\varepsilon]=\set{\lambda_0}$
(see Lemma \ref{finite}). Define the \emph{bifurcation index}
$\eta(x_0,\lambda_0)=\set{\eta_j(x_0,\lambda_0)}_{j\in\N}$
(originally defined by using topological degree for
\mbox{$\SO(2)$-equi}variant gradient mappings, see \cite{D,DR})
as follows:
\begin{equation}{\label{inddef}}
\eta_j(x_0,\lambda_0)=
\ind{\nabla H}{x_0}\cdot
\frac{\morse{Q_j((\lambda_0\!+\!\varepsilon)\nabla^2H(x_0))}\!-\!
\morse{Q_j((\lambda_0\!-\!\varepsilon)\nabla^2H(x_0))}}{2},
\end{equation}
where $\ind{\nabla H}{x_0}$ is the topological index of $x_0$ with
respect to $\nabla H,$ i.e. it is the Brouwer degree of $\nabla H$
on the neighbourhood $\Omega$ of $x_0$ such that
$\cl{\Omega}\cap\nabla H^{-1}(0)=\set{x_0}.$


Let $\Theta\equiv(0,0,\ldots).$
The following theorem gives a sufficient
condition for a trivial solution to be a branching point.


\begin{theorem}{\label{local}}
Fix isolated $x_0\in\nabla H^{-1}(0)$ and
$\lambda_0\in\Lambda(\nabla^2H(x_0)).$
If $\eta(x_0,\lambda_0)\neq\Theta$ then $(x_0,\lambda_0)$ is
a branching point of nontrivial solutions of \eqref{parham}.
Moreover, if
$U=\Omega\times(a,b)\subset \hilb\times\halfline$ is
an open bounded neighbourhood of $(x_0,\lambda_0)$ such that
$\nabla H^{-1}(0)\cap\cl{\Omega}=\set{x_0}$ and
$\Lambda(\nabla^2H(x_0))\cap[a,b]=\set{\lambda_0}$ then
$C(x_0,\lambda_0)\cap\partial U\neq\emptyset.$
\end{theorem}


The proof of Theorem \ref{local} can be obtained by a slight
modification of the proof of the following Rabinowitz-type global
bifurcation theorem for Hamiltonian systems, see \cite{DR}.


\begin{theorem}{\label{global}}
Let $\nabla H^{-1}(0)$ be finite. Fix $x_0\in\nabla H^{-1}(0),$
$\lambda_0\in\Lambda(\nabla^2H(x_0)).$
If $\eta(x_0,\lambda_0)\neq\Theta$ then $(x_0,\lambda_0)$ is
a branching point of nontrivial solutions of \eqref{parham},
and either
\begin{enumerate}
   \item $C(x_0,\lambda_0)$ is unbounded in $\hilb\times\halfline$ or
   \item $C(x_0,\lambda_0)$ is bounded and, in addition,
         $C(x_0,\lambda_0)\cap\cT(H)=\set{y_1,\ldots,y_m}$ for some
         $m\in\N,$ $y_1,\ldots,y_m\in\cT(H),$ and
         \begin{displaymath}
         \sum_{i=1}^{m}\eta(y_i)=\Theta.
         \end{displaymath}
\end{enumerate}
\end{theorem}


We will regard the set $\cS(H)=\nabla H^{-1}(0)$ of
\emph{stationary solutions} of \eqref{ham} and the set
$\cNS(H)$ of \emph{nonstationary periodic solutions}
as subsets of the Banach space $B_0\equiv B_0(\R,\R^{2n})$
of bounded functions with the supremum norm denoted as
$\norm{\cdot}_0.$ If the stationary point
$x_0\in\nabla H^{-1}(0)$ is a cluster point of $\cNS(H)$ then
it is said to be an \emph{emanation point of nonstationary periodic
solutions} of \eqref{ham}.


Denote by $\cMh\equiv\cMh(\R^{2n})$ the complete metric space of
nonempty compact subsets of $\R^{2n}$ with Hausdorff metric $d_\cH$
defined by the formula
\begin{displaymath}
d_\cH(A,B)=
\max\set{\sup_{a\in A}\dist{a}{B},\sup_{b\in B}\dist{b}{A}}
\end{displaymath}
for $A,B\in\cMh.$ For every $x_0\in\nabla H^{-1}(0),$
$x\in\cNS(H)$ we have $\norm{x-x_0}_0=d_\cH(\tr{x},x_0),$ where
$\tr{x}$ is a trajectory of $x.$ Thus $x_0$ is an emanation point of
nonstationary periodic solutions iff it is an emanation point
(in $\cMh$) of nonstationary periodic trajectories.
We say that a set $C\subset\cMh(\R^{2n})$ of nonstationary periodic
trajectories of \eqref{ham}
\emph{emanates} from $x_0$ if $x_0\in\cl{C}$ (in $\cMh(\R^{2n})$).
Note that if $C$ is connected with respect to Hausdorff metric
then the union of trajectories from $C$ is a connected subset
of $\R^{2n}.$


\begin{remarks}{\label{trajectories}}
Let $x_0$ be isolated in $\nabla H^{-1}(0).$


(1) Bifurcations of nontrivial solutions of \eqref{parham} can be
translated into emanations of nonstationary periodic trajectories
of \eqref{ham}. Namely, if $(x,\lambda)\in\hilb\times\halfline$ is
a solution of \eqref{parham} then the function
$\widetilde{x}_{\lambda}\in B_0,$
$\widetilde{x}_{\lambda}(t)=x(\frac{t}{\lambda}),$ is a solution
of \eqref{ham} with (not necessarily minimal) period $2\pi\lambda.$
(Notice that in this case $x$ can be regarded as
\mbox{$2\pi$-pe}riodic function of class $C^1$ defined on $\R.$)
Since the mapping $P:\cT(H)\cup\cNT(H)\rightarrow \cMh,$
$P(x,\lambda)=\tr{\widetilde{x}_{\lambda}}=\tr{x}=x([0,2\pi]),$
is continuous, for given connected branch
$C\subset C(x_0,\lambda_0)\cap\cNT(H)$ bifurcating from a trivial
solution $(x_0,\lambda_0)$ of \eqref{parham} we obtain a connected
branch $P(C)$ of nonstationary periodic
trajectories of \eqref{ham} emanating from $x_0.$


(2) It can be proved that there
exists a neighbourhood $U\subset\hilb\times\halfline$ of
$(x_0,\lambda_0)$ such that (not necessarily minimal) periods
$2\pi\lambda$ of trajectories
$\tr{\widetilde{x}_{\lambda}}\in P(C(x_0,\lambda_0)\cap U)$ are
arbitrarily close to $2\pi\lambda_0$ for trajectories sufficiently
close to $x_0,$ i.e. for every $\varepsilon>0$ there exists
$\delta>0$ such that if $(x,\lambda)\in C(x_0,\lambda_0)\cap U,$
$\norm{\widetilde{x}_{\lambda}-x_0}_{0}<\delta$
(or, equivalently, $d_\cH(\tr{\widetilde{x}_{\lambda}},x_0)<\delta$)
then $\abs{2\pi\lambda-2\pi\lambda_0}<\varepsilon,$ see \cite{Rd}.


(3) Moreover, if $\eta_{j}(x_0,\lambda_0)\neq 0$ then
$C(x_0,\lambda_0)\cap\cNT(H)$ contains a connected subset $C_j$
bifurcating from $(x_0,\lambda_0)$ such that for every
$(x,\lambda)\in C_j$ the number $\frac{2\pi}{j}$ is
a period (not necessarily minimal) of $x,$ which follows from the
theory of topological degree for \mbox{$\SO(2)$-equi}variant
gradient mappings. Consequently, (not necessarily minimal)
periods $\frac{2\pi\lambda}{j}$ of trajectories
$\tr{\widetilde{x}_{\lambda}}\in P(C_j\cap U)$
tend to $\frac{2\pi\lambda_0}{j}$ at $x_0.$


(4) If $\frac{\lambda_0}{\lambda}\not\in\N$
for all $\lambda\in\Lambda(\nabla^2H(x_0))\backslash\set{\lambda_0}$
then $U$ can be chosen in such a way that for every
$(x,\lambda)\in C(x_0,\lambda_0)\cap U,$ $x\neq x_0,$ the minimal
period of $\widetilde{x}_{\lambda}$ is equal to $2\pi\lambda,$
see \cite{Rd}.
\end{remarks}

\section{Algebraic results}{\label{algebraic}}


The aim of this section is to prove some algebraic lemmas which
will be used for formulation of bifurcation theorems.


Fix some $C,D\in\Sym(n,\R)$ and
$K\in\Sym(2n,\R)$ of the form
\begin{displaymath}
K=
\left[
\begin{array}{cc}
C&0\\
0&D
\end{array}
\right].
\end{displaymath}
For every $j\in\N$ define $G_j(K)\in\Sym(2n,\R)$ and
$X\in\Ort(4n,\R)$ by the formulas
\begin{displaymath}
G_j(K)=
\left[
\begin{array}{cc}
-C&jI\\
jI&-D
\end{array}
\right]
= -K+\left[
\begin{array}{cc}
0&jI\\
jI&0
\end{array}
\right],\;\;\;
X=
\left[
\begin{array}{cccc}
I&0&0&0\\
0&0&0&-I\\
0&0&I&0\\
0&I&0&0
\end{array}
\right],
\end{displaymath}
where $I\equiv I_n.$


\begin{lemma}{\label{det}}
For any $j\in\N$ we have
\begin{displaymath}
X^t Q_j(K) X=
\left[
\begin{array}{cc}
G_j(K)&0\\
0&G_j(K)
\end{array}
\right],
\end{displaymath}
$\dim\ker G_j(K)=\dim\ker[CD-j^2I]$ and $\det G_j(K)=\det[CD-j^2I].$
\end{lemma}
\begin{proof}
The first equality
of the lemma can be checked by direct calculation.
Note that the matrix
\begin{displaymath}
Y=\left[
\begin{array}{cc}
-j^{-1}I&-D\\
0&-jI
\end{array}
\right]
\end{displaymath}
is nonsingular and, moreover, $\det Y=1.$
Thus
\begin{displaymath}
\dim\ker G_j(K)=\dim\ker\left(G_j(K)Y\right)=
\dim\ker\left[
\begin{array}{cc}
j^{-1}C&CD-j^2I\\
-I&0
\end{array}
\right]
=\dim\ker[CD-j^2I]
\end{displaymath}
and, similarly,
\begin{displaymath}
\det G_j(K)=
\det\left(G_j(K)Y\right)=
\det\left[
\begin{array}{cc}
j^{-1}C&CD-j^2I\\
-I&0
\end{array}
\right]=
\det[CD-j^2I].
\end{displaymath}
\end{proof}


In what follows we assume that $A,B\in\Sym(n,\R),$
$\spect{A}=\set{\alpha_1,\ldots,\alpha_n},$
$\spect{B}=\set{\beta_1,\ldots,\beta_n},$ and
\begin{displaymath}
L=
\left[
\begin{array}{cc}
A&0\\
0&B
\end{array}
\right].
\end{displaymath}


\begin{remark}{\label{order}}
If $AB=BA$ then there exists $E\in\Ort(n,\R)$ which diagonalizes
both $A$ and $B.$ Thus we may assume without loss of generality that
\begin{displaymath}
E^tAE=J(A)=\diag{\alpha_1,\ldots,\alpha_n},\;\;\;
E^tBE=J(B)=\diag{\beta_1,\ldots,\beta_n},
\end{displaymath}
\begin{displaymath}
E^tABE=J(AB)=J(A)J(B)=\diag{\alpha_1\beta_1,\ldots,\alpha_n\beta_n}.
\end{displaymath}
We use this order of $\alpha_k,$ $\beta_k$ in the whole paper
whenever $A$ and $B$ commute.
\end{remark}


\begin{lemma}{\label{par}}
For any $j\in\N$ we have
\begin{displaymath}
\Lambda_j(L)=\set{\lambda\in\halfline\cond \det G_j(\lambda L)=0}=
\set{\frac{j}{\sqrt{\nu}}\cond \nu\in\pspect{AB}}.
\end{displaymath}
Moreover, for any fixed $\nu_0\in\pspect{AB},$ $j_0\in\N,$
the multiplicity of the root $\lambda_0=\frac{j_0}{\sqrt{\nu_0}}$ of
the polynomial $\det G_{j_0}(\lambda L)$ is equal to
$\dim\ker G_{j_0}(\lambda_0 L)=\mult{\nu_0}\leq n.$
\end{lemma}
\begin{proof}
From Lemma \ref{det} we obtain
$\det Q_j(\lambda L)=\left(\det G_j(\lambda L)\right)^2$
and
\begin{equation}{\label{roots}}
\det G_j(\lambda L)=\det[\lambda^2AB-j^2I]=
\lambda^{2n}\det\left[AB-\frac{j^2}{\lambda^2}I\right].
\end{equation}
Thus $\det G_j(\lambda L)=0$ for $\lambda\in\halfline$ if and only if
$\frac{j^2}{\lambda^2}=\nu$ for some $\nu\in\pspect{AB}.$
The equality \eqref{roots} implies also that the multiplicity of the
root $\lambda_0=\frac{j_0}{\sqrt{\nu_0}}$ is equal to $\mult{\nu_0}.$
On the other hand, by Lemma \ref{det} we have
$\dim\ker G_{j_0}(\lambda_0 L)=
\dim\ker[AB-\frac{j_0^2}{\lambda_0^2}I]=\mult{\nu_0}\leq n.$
\end{proof}


\begin{lemma}{\label{parcom}}
If $AB=BA$ then
\begin{displaymath}
\Lambda_j(L)=\set{\frac{j}{\sqrt{\alpha_k\beta_k}}
\cond\alpha_k\beta_k>0,\;k\in\set{1,\ldots,n}}
\end{displaymath}
for all $j\in\N.$
\end{lemma}
\begin{proof}
By Remark \ref{order},
$\spect{AB}=\set{\alpha_k\beta_k\cond k\in\set{1,\ldots,n}}.$
Applying Lemma \ref{par} we complete the proof.
\end{proof}


\begin{lemma}{\label{parint}}
If $\alpha_{k_0}\beta_{k_0}>0$ and
$\eigsp{A}{\alpha_{k_0}}\cap\eigsp{B}{\beta_{k_0}}\neq\set{0}$
for some $k_0\in\set{1,\ldots,n}$ then
$\frac{j}{\sqrt{\alpha_{k_0}\beta_{k_0}}}\in\Lambda_j(L)$
for every $j\in\N.$
\end{lemma}
\begin{proof}
It is easy to see that $\nu_0=\alpha_{k_0}\beta_{k_0}$ is a strictly
positive eigenvalue of $AB,$ therefore our claim is a consequence of
Lemma \ref{par}.
\end{proof}


From now on for given $\lambda_0\in\Lambda(L)$ we choose
$\varepsilon>0$ such that $\lambda_0-\varepsilon>0$ and
$\Lambda(L)\cap [\lambda_0-\varepsilon,\lambda_0+\varepsilon]=
\set{\lambda_0}$
(see Lemma \ref{finite}).


Let
\begin{equation}{\label{posneg}}
\begin{array}{c}
Y_{j}(\lambda_0)=
\set{k\in\set{1,\ldots,n}\cond \alpha_k\beta_k>0,\;
\lambda_0=\frac{j}{\sqrt{\alpha_k\beta_k}}},\\
Y_{j}^{+}(\lambda_0)=
\set{k\in Y_{j}(\lambda_0)\cond \alpha_k>0,\;\beta_k>0},\\
Y_{j}^{-}(\lambda_0)=
\set{k\in Y_{j}(\lambda_0)\cond \alpha_k<0,\;\beta_k<0}.
\end{array}
\end{equation}


For $k\in\set{1,\ldots,n}$ and $j\in\N$ define functions
$\gamma^{\pm}_{kj}:\halfline\rightarrow\R$ by the formula
\begin{displaymath}
\gamma^{\pm}_{kj}(\lambda)=
\frac{-\lambda(\alpha_k+\beta_k)
\pm\sqrt{\lambda^2(\alpha_k-\beta_k)^2+4j^2}}{2}.
\end{displaymath}


Let us gather some basic properties of $\gamma^{\pm}_{kj}.$


\begin{lemma}{\label{eigfunc}}
For every $k\in\set{1,\ldots,n},$ $j\in\N$ we have
\begin{displaymath}
\gamma^{+}_{kj}(\lambda)=0\;\wedge\;\lambda\in\halfline\;\;\;
\Leftrightarrow\;\;\;\alpha_k>0\;\wedge\;\beta_k>0\;
\wedge\;\lambda=\frac{j}{\sqrt{\alpha_k\beta_k}},
\end{displaymath}
\begin{displaymath}
\gamma^{-}_{kj}(\lambda)=0\;\wedge\;\lambda\in\halfline\;\;\;
\Leftrightarrow\;\;\;\alpha_k<0\;\wedge\;\beta_k<0\;
\wedge\;\lambda=\frac{j}{\sqrt{\alpha_k\beta_k}}.
\end{displaymath}
Moreover, for every fixed
$\lambda_0=\frac{j}{\sqrt{\alpha_k\beta_k}}$
we have
\begin{displaymath}
\alpha_k>0\;\wedge\;\beta_k>0\;\;\;\Rightarrow\;\;\;
\left\{
\begin{array}{l}
\gamma^{+}_{kj}(\lambda_0+\varepsilon)<0\;\wedge\;
\gamma^{+}_{kj}(\lambda_0-\varepsilon)>0\\
\gamma^{-}_{kj}(\lambda_0+\varepsilon)<0\;\wedge\;
\gamma^{-}_{kj}(\lambda_0-\varepsilon)<0,
\end{array}
\right.
\end{displaymath}
\begin{displaymath}
\alpha_k<0\;\wedge\;\beta_k<0\;\;\;\Rightarrow\;\;\;
\left\{
\begin{array}{l}
\gamma^{+}_{kj}(\lambda_0+\varepsilon)>0\;\wedge\;
\gamma^{+}_{kj}(\lambda_0-\varepsilon)>0\\
\gamma^{-}_{kj}(\lambda_0+\varepsilon)>0\;\wedge\;
\gamma^{-}_{kj}(\lambda_0-\varepsilon)<0.
\end{array}
\right.
\end{displaymath}
\end{lemma}
\begin{proof}
First observe that
\begin{equation}{\label{two}}
\gamma^{+}_{kj}(\lambda)\gamma^{-}_{kj}(\lambda)=
\lambda^2\alpha_k\beta_k-j^2,
\end{equation}
therefore $\gamma^{+}_{kj}(\lambda)\gamma^{-}_{kj}(\lambda)=0$ iff
$\alpha_k\beta_k>0$ and $\lambda=\lambda_0=
\frac{j}{\sqrt{\alpha_k\beta_k}}.$
If $\alpha_k>0$ and $\beta_k>0$ then $\gamma^{-}_{kj}(\lambda_0)<0,$
hence $\gamma^{+}_{kj}(\lambda_0)=0.$ For $\alpha_k<0,$ $\beta_k<0$
we have $\gamma^{+}_{kj}(\lambda_0)>0$ and thus
$\gamma^{-}_{kj}(\lambda_0)=0.$ Moreover, the derivation
of \eqref{two} gives
\begin{displaymath}
2\lambda_0\alpha_k\beta_k=
(\gamma^{+}_{kj})^\prime(\lambda_0)\gamma^{-}_{kj}(\lambda_0)+
\gamma^{+}_{kj}(\lambda_0)(\gamma^{-}_{kj})^\prime(\lambda_0),
\end{displaymath}
and thus
\begin{displaymath}
\alpha_k>0\;\;\wedge\;\;\beta_k>0\;\;\;\Rightarrow\;\;\;
(\gamma^{+}_{kj})^\prime(\lambda_0)<0,
\end{displaymath}
\begin{displaymath}
\alpha_k<0\;\;\wedge\;\;\beta_k<0\;\;\;\Rightarrow\;\;\;
(\gamma^{-}_{kj})^\prime(\lambda_0)>0.
\end{displaymath}
\end{proof}


\begin{lemma}{\label{morsecom}}
If $AB=BA$ then for every $j\in\N,$
$\lambda\in\halfline$ and fixed $\lambda_0\in\Lambda(L)$ we have
\begin{displaymath}
\spect{Q_j(\lambda L)}=\spect{G_j(\lambda L)}=
\bigcup_{k=1}^{n}
\set{\gamma^{+}_{kj}(\lambda),\gamma^{-}_{kj}(\lambda)}
\end{displaymath}
and
\begin{displaymath}
\morse{Q_j((\lambda_0+\varepsilon)L)}-
\morse{Q_j((\lambda_0-\varepsilon)L)}=
2\left(\sharp Y_{j}^{+}(\lambda_0)-\sharp Y_{j}^{-}(\lambda_0)
\right).
\end{displaymath}
\end{lemma}
\begin{proof}
Set $C=\lambda A,$ $D=\lambda B.$ In view of Lemma \ref{det},
$\spect{Q_j(\lambda L)}=\spect{G_j(\lambda L)}$ and
\begin{displaymath}
\morse{Q_j((\lambda_0+\varepsilon)L)}-
\morse{Q_j((\lambda_0-\varepsilon)L)}=
2\left(\morse{G_j((\lambda_0+\varepsilon)L)}-
\morse{G_j((\lambda_0-\varepsilon)L)}\right).
\end{displaymath}
Moreover, by Remark \ref{order} there is $E\in\Ort(n,\R)$ such that
\begin{displaymath}
E^tAE=\diag{\alpha_1,\ldots,\alpha_n},\;\;\;
E^tBE=\diag{\beta_1,\ldots,\beta_n}.
\end{displaymath}
From Lemma \ref{det} we obtain
\begin{displaymath}
\det[G_j(\lambda L)-\omega I]=\det G_j(\lambda L+\omega I)=
\det[(\lambda A+\omega I)(\lambda B+\omega I)-j^2I]=
\end{displaymath}
\begin{displaymath}
=\det[(\lambda E^tAE+\omega I)(\lambda E^tBE+\omega I)-j^2I]=
\prod_{k=1}^{n}
\left((\lambda\alpha_k+\omega)(\lambda\beta_k+\omega)-j^2\right).
\end{displaymath}
But
\begin{displaymath}
(\lambda\alpha_k+\omega)(\lambda\beta_k+\omega)-j^2=
\left(\omega-\gamma^{+}_{kj}(\lambda)\right)
\left(\omega-\gamma^{-}_{kj}(\lambda)\right)
\end{displaymath}
and thus
\begin{displaymath}
\spect{G_j(\lambda L)}=
\bigcup_{k=1}^{n}
\set{\gamma^{+}_{kj}(\lambda),\gamma^{-}_{kj}(\lambda)}.
\end{displaymath}


To compute the change of the Morse index of $G_j(\lambda L)$
observe that according to Lemma \ref{eigfunc} the eigenvalues
$\gamma^{\pm}_{kj}(\lambda)$ may change their signs at
$\lambda_0$ only if $k\in Y_j(\lambda_0),$ i.e.
$\alpha_k\beta_k>0$ and
$\lambda_0=\frac{j}{\sqrt{\alpha_k\beta_k}}.$ Using this
property, the definition of the Morse index,
the definitions of $Y_j(\lambda_0),$ $Y_{j}^{\pm}(\lambda_0),$
and the results of Lemma \ref{eigfunc} we obtain
\begin{displaymath}
\morse{G_j((\lambda_0+\varepsilon)L)}-
\morse{G_j((\lambda_0-\varepsilon)L)}=
\end{displaymath}
\begin{displaymath}
=\sharp\set{k\in\set{1,\ldots,n}\cond
\gamma^{+}_{kj}(\lambda_0+\varepsilon)<0}+
\sharp\set{k\in\set{1,\ldots,n}\cond
\gamma^{-}_{kj}(\lambda_0+\varepsilon)<0}-
\end{displaymath}
\begin{displaymath}
-\sharp\set{k\in\set{1,\ldots,n}\cond
\gamma^{+}_{kj}(\lambda_0-\varepsilon)<0}-
\sharp\set{k\in\set{1,\ldots,n}\cond
\gamma^{-}_{kj}(\lambda_0-\varepsilon)<0}=
\end{displaymath}
\begin{displaymath}
=\sharp\set{k\in Y_j(\lambda_0)\cond
\gamma^{+}_{kj}(\lambda_0+\varepsilon)<0}+
\sharp\set{k\in Y_j(\lambda_0)\cond
\gamma^{-}_{kj}(\lambda_0+\varepsilon)<0}-
\end{displaymath}
\begin{displaymath}
-\sharp\set{k\in Y_j(\lambda_0)\cond
\gamma^{+}_{kj}(\lambda_0-\varepsilon)<0}-
\sharp\set{k\in Y_j(\lambda_0)\cond
\gamma^{-}_{kj}(\lambda_0-\varepsilon)<0}=
\end{displaymath}
\begin{displaymath}
=\sharp Y_{j}^{+}(\lambda_0)+\sharp Y_{j}^{+}(\lambda_0)-
0-\sharp\left(Y_{j}^{+}(\lambda_0)\cup Y_{j}^{-}(\lambda_0)\right)=
\sharp Y_{j}^{+}(\lambda_0)-\sharp Y_{j}^{-}(\lambda_0),
\end{displaymath}
since
\begin{displaymath}
\set{k\in Y_j(\lambda_0)\cond
\gamma^{+}_{kj}(\lambda_0-\varepsilon)<0}=\emptyset
\end{displaymath}
and
\begin{displaymath}
\set{k\in Y_j(\lambda_0)\cond
\gamma^{-}_{kj}(\lambda_0-\varepsilon)<0}=
Y_{j}^{+}(\lambda_0)\cup Y_{j}^{-}(\lambda_0).
\end{displaymath}
\end{proof}


\begin{lemma}{\label{morseint}}
If $\alpha_{k_0}\beta_{k_0}>0$ and
$\eigsp{A}{\alpha_{k_0}}\cap\eigsp{B}{\beta_{k_0}}\neq\set{0}$
for some $k_0\in\set{1,\ldots,n}$ then
$\gamma^{\pm}_{k_0j}(\lambda)\in\spect{Q_j(\lambda L)}=
\spect{G_j(\lambda L)}$ for every $j\in\N,$ $\lambda\in\halfline$
and we have
\begin{displaymath}
\mult{\gamma^{\pm}_{k_0j}(\lambda)}
\geq \dim\eigsp{A}{\alpha_{k_0}}\cap\eigsp{B}{\beta_{k_0}}
\end{displaymath}
(where
$\mult{\gamma^{\pm}_{k_0j}(\lambda)}$ is the multiplicity of
$\gamma^{\pm}_{k_0j}(\lambda)$ as an eigenvalue of
$G_j(\lambda L)$).
If, in addition,
\begin{displaymath}
\dim{\eigsp{A}{\alpha_{k_0}}\cap\eigsp{B}{\beta_{k_0}}}>
\frac{1}{2}\mult{\nu_0}
\end{displaymath}
for $\nu_0=\alpha_{k_0}\beta_{k_0}\in\pspect{AB}$ then for every
fixed $j_0\in\N$ and $\lambda_0=\frac{j_0}{\sqrt{\nu_0}}$ we have
\begin{displaymath}
\morse{Q_{j_0}((\lambda_0+\varepsilon)L)}-
\morse{Q_{j_0}((\lambda_0-\varepsilon)L)}\neq 0.
\end{displaymath}
\end{lemma}
\begin{proof}
For abbreviation of notation put
$q=\dim{\eigsp{A}{\alpha_{k_0}}\cap\eigsp{B}{\beta_{k_0}}}.$
Let $E\in\Ort(n,\R)$ be such that
\begin{displaymath}
E^tAE=
\left[
\begin{array}{cc}
\alpha_{k_0}I_q&0\\
0&\widetilde{A}
\end{array}
\right],\;\;\;
E^tBE=
\left[
\begin{array}{cc}
\beta_{k_0}I_q&0\\
0&\widetilde{B}
\end{array}
\right],\;\;\;
\widetilde{A},\widetilde{B}\in\Sym(n-q,\R).
\end{displaymath}
Using Lemma \ref{det}, similarly as in the proof of
Lemma \ref{morsecom}, we obtain
\begin{displaymath}
\det[G_j(\lambda L)-\omega I]=\det G_j(\lambda L+\omega I)=
\end{displaymath}
\begin{displaymath}
=\det[(\lambda A+\omega I)(\lambda B+\omega I)-j^2I]=
\det[(\lambda E^tAE+\omega I)(\lambda E^tBE+\omega I)-j^2I]=
\end{displaymath}
\begin{displaymath}
=\left((\lambda\alpha_{k_0}+\omega)
(\lambda\beta_{k_0}+\omega)-j^2\right)^q
\cdot\det[(\lambda \widetilde{A}+\omega I)
(\lambda \widetilde{B}+\omega I)-j^2I]=
\end{displaymath}
\begin{displaymath}
=\left(\omega-\gamma^{+}_{k_0j}(\lambda)\right)^q
\left(\omega-\gamma^{-}_{k_0j}(\lambda)\right)^q
\cdot\det[(\lambda \widetilde{A}+\omega I)
(\lambda \widetilde{B}+\omega I)-j^2I]
\end{displaymath}
and thus $\gamma^{\pm}_{k_0j}(\lambda)\in\spect{G_j(\lambda L)}$
with $\mult{\gamma^{\pm}_{k_0j}(\lambda)}\geq q.$


Now, fix $j_0\in\N$ and $\lambda_0=\frac{j_0}{\sqrt{\nu_0}}.$
According to Lemma \ref{eigfunc}, if $\alpha_{k_0}>0,$
$\beta_{k_0}>0$ ($\alpha_{k_0}<0,$ $\beta_{k_0}<0$)
then the eigenvalue
$\gamma^{+}_{k_0j_0}(\lambda)$ (resp. $\gamma^{-}_{k_0j_0}(\lambda)$)
changes its sign from positive to negative (resp. from negative to
positive) when $\lambda$ changes from $\lambda_0-\varepsilon$ to
$\lambda_0+\varepsilon$ (and $\gamma^{-}_{k_0j_0}(\lambda)$
(resp. $\gamma^{+}_{k_0j_0}(\lambda)$) is nonzero).
Thus we have at least $q$ eigenvalues of $G_{j_0}(\lambda L)$
changing their signs in the same way at $\lambda_0.$ The number of
other eigenvalues changing their sings at $\lambda_0$ is less
then $q$ because $q>\frac{1}{2}\mult{\nu_0}$ and $\mult{\nu_0}$ is
the maximal number of eigenvalues of $G_{j_0}(\lambda L)$ which may
change their signs at $\lambda_0,$ since
$\dim\ker G_{j_0}(\lambda_0 L)=\mult{\nu_0},$
according to Lemma \ref{par}.
\end{proof}


\begin{lemma}{\label{morseodd}}
If $\nu_0\in\pspect{AB}$ is of odd multiplicity and
$\lambda_0=\frac{j_0}{\sqrt{\nu_0}}$ then
\begin{displaymath}
\morse{Q_{j_0}((\lambda_0+\varepsilon)L)}-
\morse{Q_{j_0}((\lambda_0-\varepsilon)L)}\neq 0.
\end{displaymath}
\end{lemma}
\begin{proof}
In view of Lemma \ref{det} it suffices to show that the product of
eigenvalues of $G_{j_0}(\lambda L)$ changes its sign at $\lambda_0.$
However, according to Lemma \ref{par}, this product is equal to
$\det G_{j_0}(\lambda L)=(\lambda-\lambda_0)^{2m+1}\psi(\lambda),$
where $2m+1$ is the multiplicity of $\nu_0$
(for some $m\in\N\cup\set{0}$) and
$\psi(\lambda_0)\neq 0,$ which completes the proof.
\end{proof}


\begin{lemma}{\label{morsepos}}
If $A$ or $B$ is strictly positive or strictly negative definite,
$\nu_0\in\pspect{AB}$ and $\lambda_0=\frac{j_0}{\sqrt{\nu_0}}$ then
for every $j\in\N$ we have
\begin{displaymath}
\morse{Q_j((\lambda_0+\varepsilon)L)}-
\morse{Q_j((\lambda_0-\varepsilon)L)}=
\end{displaymath}
\begin{displaymath}
=\left\{
\begin{array}{ll}
2\cdot s\cdot\mult{\nu}&
if\;\;\lambda_0=
\frac{j}{\sqrt{\nu}}\;\;for\;\;some\;\;\nu\in\pspect{AB},\\ \\
0&
if\;\;\lambda_0\neq
\frac{j}{\sqrt{\nu}}\;\;for\;\;all\;\;\nu\in\pspect{AB},
\end{array}
\right.
\end{displaymath}
where $s=1$ if $A$ or $B$ is strictly positive definite and
$s=-1$ if $A$ or $B$ is strictly negative definite.
\end{lemma}
\begin{proof}
Let $A$ be strictly positive or strictly negative definite and
\begin{displaymath}
M=
\left[
\begin{array}{cc}
(\sqrt{\abs{A}})^{-1}&0\\
0&\sqrt{\abs{A}}
\end{array}
\right].
\end{displaymath}
Note that $M\in\Sym(2n,\R)\cap\GL(2n,\R).$ In view of Sylvester's
law of inertia, any nonsingular transformation does not change
the Morse index of the symmetric matrix, hence
\begin{displaymath}
\morse{G_j(\lambda L)}=\morse{M^{t}G_j(\lambda L)M}=
\morse{MG_j(\lambda L)M}=
\end{displaymath}
\begin{displaymath}
=\morse{
\left[
\begin{array}{cc}
-\lambda\,\sgn(A)I&jI\\
jI&-\lambda\sqrt{\abs{A}}B\sqrt{\abs{A}}
\end{array}
\right]}.
\end{displaymath}
The matrix $\sqrt{\abs{A}}B\sqrt{\abs{A}}$ is symmetric, so it has
$n$ real eigenvalues $\omega_1,\ldots,\omega_n.$ On the other hand,
these eigenvalues are exactly those of $\abs{A}B$ because
\begin{displaymath}
\det[\sqrt{\abs{A}}B\sqrt{\abs{A}}-\omega I]
=\det\abs{A}\cdot\det[B-\omega\abs{A}^{-1}]=
\det[\abs{A}B-\omega I].
\end{displaymath}
Applying Lemma \ref{morsecom} with $\sgn(A)I$ and
$\sqrt{\abs{A}}B\sqrt{\abs{A}}$ instead of $A$ and $B,$ respectively,
we get
\begin{displaymath}
\morse{G_j((\lambda_0+\varepsilon)L)}-
\morse{G_j((\lambda_0-\varepsilon)L)}=
\end{displaymath}
\begin{displaymath}
=\sgn(A)\cdot\sharp\set{k\in\set{1,\ldots,n}\cond
\sgn(A)\omega_k>0\;\wedge\;
\lambda_0=\frac{j}{\sqrt{\sgn(A)\omega_k}}}.
\end{displaymath}
But $\sgn(A)\omega_1,\ldots,\sgn(A)\omega_n$ are the eigenvalues of
$\sgn(A)\abs{A}B=AB,$ which completes the proof for $A$ strictly
positive or negative. If $B$ is strictly positive or negative
definite, consider
\begin{displaymath}
M=
\left[
\begin{array}{cc}
\sqrt{\abs{B}}&0\\
0&(\sqrt{\abs{B}})^{-1}
\end{array}
\right].
\end{displaymath}
\end{proof}


Note that the number $s$ in the above lemma is well defined because
we assume that $\pspect{AB}\neq\emptyset.$ If $A$
and $B$ were nonsingular and of different signs then $AB$ would be
strictly negative definite, a contradiction. For example, if $A$ is
strictly positive and $B$ is strictly negative then
$\spect{AB}=\spect{\sqrt{A}B\sqrt{A}}$ and for every $v\in\R^n$
we have
\begin{displaymath}
\left\langle \sqrt{A}B\sqrt{A}v,v \right\rangle=
\left\langle B\sqrt{A}v,\sqrt{A}v \right\rangle <0,
\end{displaymath}
where $\langle\cdot,\cdot\rangle$ is an inner product in $\R^n.$


\section{Local bifurcations}


In this section we formulate local bifurcation theorems for
autonomous Hamiltonian systems (with block-diagonal Hessian
of the Hamiltonian at a stationary point)
in terms of the topological degree
of the gradient of the Hamiltonian and eigenvalues
of its Hessian computed at a stationary point.


Assume that $H\in C^2(\R^n\times\R^n,\R)$ and
that for fixed $x_0\in\nabla H^{-1}(0),$
isolated in $\nabla H^{-1}(0),$ the Hessian of $H$ at $x_0$
has the form
\begin{equation}{\label{localform}}
\nabla^2H(x_0)=
\left[
\begin{array}{cc}
A&0\\
0&B
\end{array}
\right]
\end{equation}
for some $A,B\in\Sym(n,\R).$


Let $\spect{A}=\set{\alpha_1,\ldots,\alpha_n},$
$\spect{B}=\set{\beta_1,\ldots,\beta_n}$ and assume the
convention of Remark \ref{order} for the order of
$\alpha_k,$ $\beta_k.$


In view of Theorem \ref{nec} we may suspect
that the point $(x_0,\lambda_0)$ is a bifurcation point of
nontrivial solutions of \eqref{parham} provided that
$\lambda_0\in\Lambda(\nabla^2H(x_0)).$ By Lemma \ref{par}
this means that
\begin{equation}{\label{explnec}}
\lambda_0=
\frac{j_0}{\sqrt{\nu_0}},\;\;j_0\in\N,\;\;
\nu_0\in\pspect{AB}.
\end{equation}
If $AB=BA$ then every such a $\lambda_0$ can be written
as $\lambda_0=\frac{j_0}{\sqrt{\alpha_{k_0}\beta_{k_0}}}$ for some
$k_0\in\set{1,\ldots,n}$ such that $\alpha_{k_0}\beta_{k_0}>0$
(see Lemma \ref{parcom}). In the case of $AB\neq BA$ we cannot
write $\lambda_0$ in this form in general, but it is possible
to do it if $\alpha_{k_0}\beta_{k_0}>0$ and
$V_A(\alpha_{k_0})\cap V_B(\beta_{k_0})\neq\set{0}$
(see Lemma \ref{parint}).


The following conditions will be used in theorems of this section.


\noindent (A1) $H\in C^2(\R^n\times\R^n,\R),$
$x_0\in\nabla H^{-1}(0)$ is isolated in $\nabla H^{-1}(0)$ and
$\ind{\nabla H}{x_0}\neq 0,$


\noindent (A2) $\displaystyle
\nabla^2H(x_0)=
\left[
\begin{array}{cc}
A&0\\
0&B
\end{array}
\right],\;\;\;A,B\in\Sym(n,\R),\;\;\;$
$\nu_0\in\pspect{AB},$


\noindent (A3) $\lambda_0=\frac{j_0}{\sqrt{\nu_0}},$ $j_0\in\N,$


\noindent (A4) $AB=BA$ and
$\sharp Y_{j}^{+}(\lambda_0)\neq\sharp Y_{j}^{-}(\lambda_0)$ for some
$j\in\N,$ where $Y_{j}^{+}(\lambda_0)$ and $Y_{j}^{-}(\lambda_0)$
are given by \eqref{posneg}.


\begin{theorem}{\label{localcom}}
Assume that all conditions (A1)-(A4) are satisfied.
Then $(x_0,\lambda_0)$ is a branching point of nontrivial
solutions of
\begin{displaymath}
\left\{
\begin{array}{l}
\dot{x}(t)=\lambda J\nabla H(x(t))\\
x(0)=x(2\pi),
\end{array}
\right.
\end{displaymath}
where $\lambda\in\halfline.$
Moreover, if
$U=\Omega\times(a,b)\subset\hilb\times\halfline$ is an open bounded
neighbourhood of $(x_0,\lambda_0)$ such that
$\nabla H^{-1}(0)\cap\cl{\Omega}=\set{x_0}$ and $(x_0,\lambda_0)$ is
the only trivial solution in $\cl{U}$ satisfying \eqref{explnec}
then $C(x_0,\lambda_0)\cap\partial U\neq\emptyset.$
\end{theorem}
\begin{proof}
Observe that $\lambda_0\in\Lambda(\nabla^2H(x_0))$ (see (A2), (A3),
and \eqref{explnec}). Moreover, the assumptions of
Lemma \ref{morsecom} are satisfied for $L=\nabla^2H(x_0),$ according
to (A2) and (A4). Using this lemma and equality
\eqref{inddef} we obtain
$\eta_j(x_0,\lambda_0)=
\ind{\nabla H}{x_0}\cdot
\left(\sharp Y_{j}^{+}(\lambda_0)-\sharp Y_{j}^{-}(\lambda_0)\right)$
and thus $\eta_j(x_0,\lambda_0)\neq 0,$ by (A1) and (A4).
Applying Theorem \ref{local} we complete the proof.
\end{proof}


The corresponding theorem concerning emanations of periodic
trajectories from a stationary point can be formulated as follows.


\begin{theorem}{\label{emancom}}
Assume that all conditions (A1)-(A4) are satisfied.
Then there exists a connected (with respect to
Hausdorff metric) set of nonstationary periodic trajectories of
\begin{displaymath}
\dot{x}(t)=J\nabla H(x(t))
\end{displaymath}
emanating from
$x_0$ with (not necessarily minimal) periods arbitrarily close to
$\frac{j_0}{j}\frac{2\pi}{\sqrt{\nu_0}}$ for
trajectories sufficiently close to $x_0.$
Moreover, if $j_0=1$ and
$\sqrt{\frac{\nu}{\nu_0}}\not\in\N$ for all
$\nu\in\pspect{AB}\backslash\set{\nu_0}$ then
there exists a connected set of nonstationary periodic
trajectories emanating from
$x_0$ with minimal periods arbitrarily close to
$\frac{2\pi}{\sqrt{\nu_0}}$ for trajectories
sufficiently close to $x_0.$
\end{theorem}
\begin{proof}
Our claim is a consequence of Theorem \ref{localcom},
Lemma \ref{morsecom}, and Remarks \ref{trajectories}.
To obtain minimal periods of trajectories observe that if
$j_0=1$ and $\sqrt{\frac{\nu}{\nu_0}}\not\in\N$
for all $\nu\in\pspect{AB}\backslash\set{\nu_0}$ then
$\frac{\lambda_0}{\lambda}\not\in\N$ for all
$\lambda\in\Lambda(\nabla^2H(x_0))\backslash\set{\lambda_0}.$
\end{proof}


Notice that if $j_0=1$ and $\sqrt{\frac{\nu}{\nu_0}}\not\in\N$ for
all $\nu\in\pspect{AB}\backslash\set{\nu_0}$ then condition (A4)
in the above theorem can be satisfied only for $j=1,$ since
in this case $\lambda_0=\frac{1}{\sqrt{\nu_0}}=
\frac{j}{\sqrt{\alpha_k\beta_k}}$ only
for $j=j_0=1$ and $\nu_0=\alpha_k\beta_k,$ see \eqref{posneg}.


The following theorem can be proved in the same way as
Theorem \ref{localcom}, by using
Lemma \ref{morseint} instead of \ref{morsecom}.


\begin{theorem}{\label{localeigen}}
Suppose that conditions (A1), (A2), (A3) are
satisfied and that for some $k_0\in\set{1,\ldots,n}$ we have
\begin{equation}{\label{eigencond}}
\nu_0=\alpha_{k_0}\beta_{k_0},\;\;
\dim V_A(\alpha_{k_0})\cap V_B(\beta_{k_0})>\frac{1}{2}\mult{\nu_0}.
\end{equation}
Then the conclusion of Theorem \ref{localcom} holds true.
\end{theorem}


Combining Theorem \ref{localeigen} and Lemma \ref{morseint}
(for $j_0=1$) with Remarks \ref{trajectories} we obtain
the corresponding emanation result.


\begin{theorem}{\label{emaneigen}}
Suppose that assumptions (A1), (A2) are satisfied and that
condition \eqref{eigencond} is fulfilled for some
$k_0\in\set{1,\ldots,n}.$
Then there exists a connected (with respect to
Hausdorff metric) set of nonstationary periodic trajectories of
\begin{displaymath}
\dot{x}(t)=J\nabla H(x(t))
\end{displaymath}
emanating from
$x_0$ with (not necessarily minimal) periods arbitrarily close to
$\frac{2\pi}{\sqrt{\nu_0}}$ for
trajectories sufficiently close to $x_0.$
Moreover, if $\sqrt{\frac{\nu}{\nu_0}}\not\in\N$ for all
$\nu\in\pspect{AB}\backslash\set{\nu_0}$ then
there exists a connected set of nonstationary periodic
trajectories emanating from
$x_0$ with minimal periods arbitrarily close to
$\frac{2\pi}{\sqrt{\nu_0}}$ for trajectories
sufficiently close to $x_0.$
\end{theorem}


Similarly as above, application of Lemma \ref{morseodd} gives us
the following.


\begin{theorem}{\label{localodd}}
Assume that conditions (A1), (A2) are
satisfied and the multiplicity of $\nu_0$ is odd.
Then
\begin{enumerate}
   \item the conclusion of Theorem \ref{localcom} is true
   for every $\lambda_0$ satisfying (A3),
   \item the conclusion of Theorem \ref{emaneigen} holds.
\end{enumerate}
\end{theorem}


Finally, applying Lemma \ref{morsepos} (for $j=j_0$) we obtain


\begin{theorem}{\label{localpos}}
Let conditions (A1) and (A2) be fulfilled.
Suppose that $A$ or $B$ is strictly positive or strictly negative
definite. Then
\begin{enumerate}
   \item the conclusion of Theorem \ref{localcom} is true
   for every $\lambda_0$ satisfying (A3),
   \item the conclusion of Theorem \ref{emaneigen} holds.
\end{enumerate}
\end{theorem}


\begin{corollary}{\label{localmin}}
Let $H\in C^2(\R^n\times\R^n,\R)$ admit a strict local minimum or
maximum at $x_0$ and
\begin{displaymath}
\nabla^2H(x_0)=
\left[
\begin{array}{cc}
A&0\\
0&B
\end{array}
\right],\;\;\;A,B\in\Sym(n,\R).
\end{displaymath}
If $A$ is nonsingular and $B\neq 0,$ or
$B$ is nonsingular and $A\neq 0,$ then
\begin{enumerate}
   \item $\pspect{AB}\neq\emptyset,$
   \item the conclusion of Theorem \ref{localcom} holds for any
   $\lambda_0=\frac{j_0}{\sqrt{\nu_0}},$ $j_0\in\N,$
   $\nu_0\in\pspect{AB},$
   \item the conclusion of Theorem \ref{emaneigen}
   is true for every $\nu_0\in\pspect{AB}.$
\end{enumerate}
\end{corollary}
\begin{proof}
Since $H$ admits a strict local minimum (maximum) at $x_0,$ we have
$x_0\in \nabla H^{-1}(0),$ $x_0$ is isolated in $\nabla H^{-1}(0),$
$\ind{\nabla H}{x_0}=1\neq 0$
(see \cite{A})  and $\nabla^2H(x_0)$
is nonnegative (resp. nonpositive) definite, hence $\sgn(A)=\sgn(B).$
If, for example, $A$ is nonsingular then
\begin{displaymath}
\pspect{AB}=\pspect{\abs{A}\abs{B}}=
\pspect{\sqrt{\abs{A}}\abs{B}\sqrt{\abs{A}}}=
\pspect{\sqrt{\abs{A}}^t\abs{B}\sqrt{\abs{A}}}.
\end{displaymath}
But $B\neq 0,$ therefore
$\pspect{AB}\neq 0,$ in view of Sylvester's law of inertia. Applying
Theorem \ref{localpos} we complete the proof.
\end{proof}


\begin{example}
Consider $H:\R^3\times\R^3\rightarrow\R$ given by the formula
\begin{displaymath}
H(x)=H(x_1,\ldots,x_6)=
x_1^2+2(x_2-1)^2-x_3^2+(x_4+(x_3-2)^2)^6+x_5^2+(x_6-1)^2-x_5x_6.
\end{displaymath}
In this case we have
$\nabla H^{-1}(0)=\set{(0,1,0,-4,\frac{2}{3},\frac{4}{3})}$
and $\ind{\nabla H}{(0,1,0,-4,\frac{2}{3},\frac{4}{3})}=-1.$ (The
last equality can be obtained by using an algorithm described
in \cite{LS}.)
The Hessian $\nabla^2H(0,1,0,-4,\frac{2}{3},\frac{4}{3})$ has
the block-diagonal form \eqref{localform} with
\begin{displaymath}
A=
\left[
\begin{array}{rrr}
2&0&0\\
0&4&0\\
0&0&-2
\end{array}
\right],\;\;\;
B=
\left[
\begin{array}{rrr}
0&0&0\\
0&2&-1\\
0&-1&2
\end{array}
\right].
\end{displaymath}
Moreover, $AB\neq BA$ and
$\spect{AB}=\set{0,2+2\sqrt{7},2-2\sqrt{7}}.$
According to Theorem \ref{localodd}, for every $j_0\in\N$ the point
$\left((0,1,0,-4,\frac{2}{3},\frac{4}{3}),
\frac{j_0}{\sqrt{2+2\sqrt{7}}}\right)$ is
a branching point of nontrivial solutions of \eqref{parham} and there
exists a connected set of nonstationary periodic trajectories
of \eqref{ham} emanating from $(0,1,0,-4,\frac{2}{3},\frac{4}{3})$
with minimal periods tending to $\frac{2\pi}{\sqrt{2+2\sqrt{7}}}$ at
$(0,1,0,-4,\frac{2}{3},\frac{4}{3}).$
\end{example}


\begin{example}{\label{exmin}}
Define $H:\R^2\times\R^2\rightarrow\R$ as
\begin{displaymath}
H(x)=H(x_1,\ldots,x_4)=x_1^2+2x_2^2+x_3^2+(x_4+(x_3-1)^3)^4.
\end{displaymath}
Observe that $\nabla H^{-1}(0)=\set{(0,0,0,1)}$ and
$H$ admits a strict minimum at $(0,0,0,1).$ Moreover,
$\nabla^2H(0,0,0,1)$ is of the form \eqref{localform} with
\begin{displaymath}
A=
\left[
\begin{array}{cc}
2&0\\
0&4
\end{array}
\right],\;\;\;
B=
\left[
\begin{array}{cc}
2&0\\
0&0
\end{array}
\right],\;\;\;
AB=BA=
\left[
\begin{array}{cc}
4&0\\
0&0
\end{array}
\right].
\end{displaymath}
By Corollary \ref{localmin}, for every $j_0\in\N$ the point
$\left((0,0,0,1),
\frac{j_0}{2}\right)$ is
a branching point of nontrivial solutions of \eqref{parham} and there
exists a connected set of nonstationary periodic trajectories
of \eqref{ham} emanating from $(0,0,0,1)$ with
minimal periods tending to $\pi$ at
$(0,0,0,1).$
\end{example}


\section{Global bifurcations}


As in the previous section suppose that
$H\in C^2(\R^n\times\R^n,\R),$ but $\nabla H^{-1}(0)$ is finite.
Assume also that for every fixed $x_0\in\nabla H^{-1}(0)$ we have
\begin{equation}{\label{globalform}}
\nabla^2H(x_0)=L(x_0)=
\left[
\begin{array}{cc}
A(x_0)&0\\
0&B(x_0)
\end{array}
\right],\;\;\;A(x_0),B(x_0)\in\Sym(n,\R).
\end{equation}
According to Lemma \ref{par}, the set of parameters $\lambda_0$
which satisfy the bifurcation necessary condition from Theorem
\ref{nec} is equal to $\displaystyle\Lambda(\nabla^2H(x_0))=
\bigcup_{j\in\N}\Lambda_j(\nabla^2H(x_0)),$ where
\begin{displaymath}
\Lambda_j(\nabla^2H(x_0))=
\set{\frac{j}{\sqrt{\nu}}\cond \nu\in\pspect{A(x_0)B(x_0)}}.
\end{displaymath}
Let
\begin{displaymath}
\spect{A(x_0)}=\set{\alpha_1(x_0),\ldots,\alpha_n(x_0)},\;\;\;
\spect{B(x_0)}=\set{\beta_1(x_0),\ldots,\beta_n(x_0)}.
\end{displaymath}
Repeat Remark \ref{order} and equalities \eqref{posneg} for
$A\equiv A(x_0),$ $B\equiv B(x_0),$ $\alpha_k\equiv\alpha_k(x_0),$
$\beta_k\equiv\beta_k(x_0)$ and
$Y_{j}^{\pm}(\lambda_0)\equiv Y_{j}^{\pm}(x_0,\lambda_0).$


\begin{theorem}{\label{globalcom}}
Assume that $H\in C^2(\R^n\times\R^n,\R),$ $\nabla H^{-1}(0)$ is
finite, and that for every $\xi\in\nabla H^{-1}(0)$ we have
\begin{displaymath}
\nabla^2H(\xi)=
\left[
\begin{array}{cc}
A(\xi)&0\\
0&B(\xi)
\end{array}
\right],\;\;\;A(\xi),B(\xi)\in\Sym(n,\R),\;\;\;
A(\xi)B(\xi)=B(\xi)A(\xi).
\end{displaymath}
Let $x_0\in\nabla H^{-1}(0)$ and $\ind{\nabla H}{x_0}\neq 0.$
If $\lambda_0=\frac{j_0}{\sqrt{\nu_0}},$ $j_0\in\N,$
$\nu_0\in\pspect{A(x_0)B(x_0)},$ and
$\sharp Y_{l}^{+}(x_0,\lambda_0)\neq
\sharp Y_{l}^{-}(x_0,\lambda_0)$ for some $l\in\N$
then $(x_0,\lambda_0)$
is a branching point of nontrivial solutions of
\begin{displaymath}
\left\{
\begin{array}{l}
\dot{x}(t)=\lambda J\nabla H(x(t))\\
x(0)=x(2\pi),
\end{array}
\right.
\end{displaymath}
where $\lambda\in\halfline,$ and either
\begin{enumerate}
   \item $C(x_0,\lambda_0)$ is unbounded in $\hilb\times\halfline$ or
   \item $C(x_0,\lambda_0)$ is bounded and, in addition,
         $C(x_0,\lambda_0)\cap\cT(H)=
         \set{(\xi_1,\eta_1),\ldots,(\xi_m,\eta_m)}$
         for some $m\in\N,$ $\xi_i\in\nabla H^{-1}(0),$
         $\eta_i=\frac{j_i}{\sqrt{\omega_i}},$ $j_i\in\N,$
         $\omega_i\in\pspect{A(\xi_i)B(\xi_i)},$ $i=1,\ldots,m,$
         and for any $j\in\N$ we have
         \begin{displaymath}
         \sum_{i=1}^{m}\ind{\nabla H}{\xi_i}\cdot
         (\sharp Y_{j}^{+}(\xi_i,\eta_i)-
         \sharp Y_{j}^{-}(\xi_i,\eta_i))=0.
         \end{displaymath}
\end{enumerate}
\end{theorem}
\begin{proof}
According to Theorem \ref{global}, the sum of bifurcation indices of
the points from $C(x_0,\lambda_0)\cap\cT(H)$ is equal to $\Theta.$
But
$\eta_j(\xi_i,\eta_i)=
\ind{\nabla H}{\xi_i}\cdot
(\sharp Y_{j}^{+}(\xi_i,\eta_i)-
\sharp Y_{j}^{-}(\xi_i,\eta_i))$ for $i=1,\ldots,m,$
in view of equality \eqref{inddef} and Lemma \ref{morsecom}.
\end{proof}


\begin{theorem}{\label{globalpos}}
Assume that $H\in C^2(\R^n\times\R^n,\R),$ $\nabla H^{-1}(0)$ is
finite, and that for every $\xi\in\nabla H^{-1}(0)$ we have
\begin{displaymath}
\nabla^2H(\xi)=
\left[
\begin{array}{cc}
A(\xi)&0\\
0&B(\xi)
\end{array}
\right],\;\;\;A(\xi),B(\xi)\in\Sym(n,\R),
\end{displaymath}
where $A(\xi)$ or $B(\xi)$ is strictly positive or strictly
negative definite.
Let $x_0\in\nabla H^{-1}(0)$ and $\ind{\nabla H}{x_0}\neq 0.$
If $\lambda_0=\frac{j_0}{\sqrt{\nu_0}},$ $j_0\in\N,$
$\nu_0\in\pspect{A(x_0)B(x_0)}$ then $(x_0,\lambda_0)$
is a branching point of nontrivial solutions of
\begin{displaymath}
\left\{
\begin{array}{l}
\dot{x}(t)=\lambda J\nabla H(x(t))\\
x(0)=x(2\pi),
\end{array}
\right.
\end{displaymath}
where $\lambda\in\halfline,$ and either
\begin{enumerate}
   \item $C(x_0,\lambda_0)$ is unbounded in $\hilb\times\halfline$ or
   \item $C(x_0,\lambda_0)$ is bounded and, in addition,
         $C(x_0,\lambda_0)\cap\cT(H)$ is finite, for any $j\in\N$ we
         have
         \begin{displaymath}
         C(x_0,\lambda_0)\cap\bigcup_{\xi\in\nabla H^{-1}(0)}
         \set{\xi}\times\Lambda_j(\nabla^2H(\xi))=
         \set{\left(\xi_1,\frac{j}{\sqrt{\omega_1}}\right),\ldots,
         \left(\xi_m,\frac{j}{\sqrt{\omega_m}}\right)}
         \end{displaymath}
         for some $m\in\N,$ $\xi_i\in\nabla H^{-1}(0),$
         $\omega_i\in\pspect{A(\xi_i)B(\xi_i)},$ $i=1,\ldots,m$
         (if the above intersection is nonempty),
         and
         \begin{displaymath}
         \sum_{i=1}^{m}\ind{\nabla H}{\xi_i}
         \cdot s(\xi_i)\cdot \mult{\omega_i}=0,
         \end{displaymath}
         where $s(\xi_i)=1$ if $A(\xi_i)$ or $B(\xi_i)$ is strictly
         positive definite and $s(\xi_i)=-1$ if $A(\xi_i)$ or
         $B(\xi_i)$ is strictly negative definite.
\end{enumerate}
\end{theorem}
\begin{proof}
According to Theorem \ref{global}, the sum of bifurcation indices of
the points from $C(x_0,\lambda_0)\cap\cT(H)$ is equal to $\Theta,$
i.e. it vanishes at every coordinate.
In view of the equality \eqref{inddef} and Lemma \ref{morsepos} the
set
\begin{displaymath}
C(x_0,\lambda_0)\cap\bigcup_{\xi\in\nabla H^{-1}(0)}\set{\xi}\times
\Lambda_j(\nabla^2H(\xi))
\end{displaymath}
consists of those points from $C(x_0,\lambda_0)\cap\cT(H)$
for which the $j$th coordinate of the bifurcation index
$\eta$ can be nonzero. Namely,
$\eta_j(\xi_i,\frac{j}{\sqrt{\omega_i}})=
\ind{\nabla H}{\xi_i}
\cdot s(\xi_i)\cdot\mult{\omega_i}$ for $i=1,\ldots,m.$
\end{proof}


\begin{corollary}{\label{globalone}}
If the assumptions of Theorem \ref{globalpos} are satisfied and
$\nabla H^{-1}(0)=\set{x_0}$ then $C(x_0,\lambda_0)$ is unbounded
in $\hilb\times\halfline.$
\end{corollary}


\begin{example}
Let $H$ be such as in Example \ref{exmin}. Then $A$ is strictly
positive definite and $\ind{\nabla H}{(0,0,0,1)}=1.$ In view of
Corollary \ref{globalone}, for every $j_0\in\N$ the connected
branch $C\left((0,0,0,1),\frac{j_0}{2}\right)$ bifurcating from
$\left((0,0,0,1),\frac{j_0}{2}\right)$ is unbounded in
$\hilb\times\halfline.$
\end{example}


From now on we assume that $H\in C^2(\R^n\times\R^n,\R)$
satisfies the assumptions of Theorem \ref{globalpos}.
Let $s(\xi)$ be such as in that theorem and set
\begin{displaymath}
S_{+}(H)=\set{\xi\in\nabla H^{-1}(0)\cond
\ind{\nabla H}{\xi}\cdot s(\xi)>0},
\end{displaymath}
\begin{displaymath}
S_{-}(H) =\set{\xi\in\nabla H^{-1}(0)\cond
\ind{\nabla H}{\xi}\cdot s(\xi)<0},
\end{displaymath}
\begin{displaymath}
p(H)=\set{(\xi,\omega)\cond
\xi\in S_{+}(H),\omega\in\pspect{A(\xi)B(\xi)}},
\end{displaymath}
\begin{displaymath}
n(H)=\set{(\xi,\omega)\cond
\xi\in S_{-}(H), \omega\in\pspect{A(\xi)B(\xi)}},
\end{displaymath}
\begin{displaymath}
\cE(H)=\sum_{(\xi,\omega)\in p(H)\cup n(H)}
\ind{\nabla H}{\xi}\cdot s(\xi)\cdot\mult{\omega}.
\end{displaymath}


Let us formulate further corollaries to
Theorem \ref{globalpos}. If $A(\xi)=I$ for all
$\xi\in\nabla H^{-1}(0)$ then they
imply corresponding corollaries from \cite{Rd}.


\begin{corollary}
If $\cE(H)\neq 0$ then for every $j\in\N$ there exists
$(\xi,\omega)\in p(H)\cup n(H)$ such that
$C(\xi,\frac{j}{\sqrt{\omega}})$ is unbounded in
$\hilb\times\halfline.$ Moreover, if additionally $p(H)=\emptyset$
or $n(H)=\emptyset$ then
$C(\xi,\frac{j}{\sqrt{\omega}})$ are unbounded for all $j\in\N,$
$(\xi,\omega)\in p(H)\cup n(H).$
\end{corollary}
\begin{proof}
Fix $j\in\N$ and observe that if for all
$(\xi,\omega)\in p(H)\cup n(H)$ sets $C(\xi,\frac{j}{\sqrt{\omega}})$
were bounded then the sum of
$\ind{\nabla H}{\xi}\cdot s(\xi)\cdot\mult{\omega}$
over $p(H)\cup n(H)$ would be equal to $0,$ a contradiction.
\end{proof}


\begin{corollary}
If $\cE(H)\neq 0$ and
$\abs{\ind{\nabla H}{\xi}\cdot\mult{\omega}}=c=const$ for all
$(\xi,\omega)\in p(H)\cup n(H)$ then for every $j\in\N$ sets
$C(\xi,\frac{j}{\sqrt{\omega}})$ are unbounded in
$\hilb\times\halfline$ for at least $\abs{\sharp p(H)-\sharp n(H)}$
of $(\xi,\omega)\in p(H)\cup n(H).$
\end{corollary}
\begin{proof}
Assume, for example, that $\sharp p(H)>\sharp n(H).$
Denote by $Z_p\;$ ($Z_n$) the set of such points
$(\xi,\omega)\in p(H)$ (resp. $(\xi,\omega)\in n(H)$) that
$C(\xi,\frac{j}{\sqrt{\omega}})$ is bounded. The sum of
$\ind{\nabla H}{\xi}\cdot s(\xi)\cdot\mult{\omega}$ over
$Z_p\cup Z_n$ is equal to $0,$ in view of Theorem \ref{globalpos}.
Thus $\sharp Z_p=\sharp Z_n.$ But $\sharp Z_n\leq\sharp n(H),$ hence
the number of $(\xi,\omega)\in p(H)\cup n(H)$ for which
$C(\xi,\frac{j}{\sqrt{\omega}})$ is unbounded is equal to
$(\sharp p(H)-\sharp Z_p)+(\sharp n(H)-\sharp Z_n)=
\sharp p(H)+\sharp n(H)-2\sharp Z_n\geq
\sharp p(H)+\sharp n(H)-2\sharp n(H)=\sharp p(H)-\sharp n(H).$
\end{proof}


Obviously, unbounded sets $C(\xi,\frac{j}{\sqrt{\omega}})$ from the
above corollary need not be different for different
$(\xi,\omega)\in p(H)\cup n(H).$


\begin{corollary}
Suppose that $\deg(\nabla H,U,0)\neq 0$ for some open and bounded
$U\subset\R^{2n}$ such that $\nabla H^{-1}(0)\subset U.$
Let $p(H)\cup n(H)\neq\emptyset.$ If
$\sharp\pspect{\nabla^2H(\xi)}=b=const,$ $s(\xi)=s=const$
and $\mult{\omega}=m=const$
for all $(\xi,\omega)\in p(H)\cup n(H)$ then for every $j\in\N$
there exists $(\xi,\omega)\in p(H)\cup n(H)$ such that
$C(\xi,\frac{j}{\sqrt{\omega}})$ is unbounded in
$\hilb\times\halfline.$
\end{corollary}
\begin{proof}
Suppose, contrary to our claim, that
for some $j\in\N$ sets $C(\xi,\frac{j}{\sqrt{\omega}})$ are
bounded for all $(\xi,\omega)\in p(H)\cup n(H).$ According to
Theorem \ref{globalpos} the sum of
$\ind{\nabla H}{\xi}\cdot s(\xi)\cdot\mult{\omega}$ over these points
is equal to $0.$ On the other hand this sum is equal to
\begin{displaymath}
b\cdot s\cdot m\cdot\sum_{\xi\in\nabla H^{-1}(0)}\ind{\nabla H}{\xi}=
b\cdot s\cdot m\cdot\deg(\nabla H,U,0)\neq 0,
\end{displaymath}
a contradiction.
\end{proof}


Note that the condition $\deg(\nabla H,U,0)\neq 0$ in the above
corollary is satisfied if $H(x)\rightarrow +\infty$ as
$\abs{x}\rightarrow +\infty$ (see \cite{A}). For strictly
convex $H$ the last condition is equivalent
to the condition $\nabla H^{-1}(0)\neq\emptyset$ (see \cite{MW}).



\end{document}